\documentclass[11pt]{article}
\usepackage{epsfig,makeidx,amsmath,amsfonts,latexsym,subfigure}
\newtheorem{theorem}{Theorem}[section]
\newtheorem{prop}{Proposition}[section]
\newtheorem{lemma}{Lemma}[section]
\newtheorem{conj}{Conjecture}[section]

\begin{document}

\title{\LARGE\textbf{{Coexistence in a two-type continuum growth model}}}

\author{Maria Deijfen \thanks{Stockholm University. E-mail: mia@matematik.su.se}
\and Olle H\"{a}ggstr\"{o}m \thanks{Chalmers University of
Technology. E-mail: olleh@math.chalmers.se}}

\date{August 2004}

\maketitle

\thispagestyle{empty}

\begin{abstract}

\noindent We consider a stochastic model, describing the growth of
two competing infections on $\mathbb{R}^d$. The growth takes place
by way of spherical outbursts in the infected region, an outburst
in the type 1 (2) infected region causing all previously
uninfected points within a stochastic distance from the outburst
location to be type 1 (2) infected. The main result is that, if
the infection types have the same intensity, then there is a
strictly positive probability that both infection types grow
unboundedly.

\vspace{1cm}

\noindent \emph{Keywords:} Richardson's model,
competing growth, coexistence.

\vspace{0.5cm}

\noindent AMS 2000 Subject Classification: Primary 60K35\newline
\hspace*{5.15cm} Secondary 82B43.
\end{abstract}

\section{Introduction}

\noindent In Deijfen (2003), a growth model is introduced that
describes the spread of an entity, henceforth referred to as an
infection, on $\mathbb{R}^d$. The model is generalized in Deijfen
\emph{et al} (2003) to encompass two infection types, making it a
model for competition on $\mathbb{R}^d$. Deijfen \emph{et al}
(2003) is mainly concerned with the case when the infections have
unequal intensities and the main result is roughly that two
infections with different intensities cannot simultaneously grow
to occupy infinite parts of $\mathbb{R}^d$. In this paper we
consider infections with the same intensity and show that in this
case there is a positive probability that both infection types
reach points arbitrarily far from the origin.\medskip

\noindent The Richardson model, introduced in Richardson (1973),
describes growth on $\mathbb{Z}^d$. Sites are either healthy or
infected, the dynamics being that a healthy site is infected at a
rate proportional to the number of infected neighbors and once
infected it remains so forever. The model of Deijfen (2003) can be
viewed as a kind of continuum analog of the Richardson model.
Starting from an arbitrary initial set, the growth is driven by
randomly sized spherical outbursts in the infected region, whose
value at time $t$ is denoted by $S_t$. The time until an outburst
takes place is exponentially distributed with parameter
proportional to $|S_t|$, where $|\cdot|$ denotes Lebesgue measure,
and the outburst point is uniformly distributed over $S_t$. When
an outburst occurs it causes the previously uninfected parts of a
ball centered at the outburst point to be infected, the radii in
the outburst balls being i.i.d.\ random variables with
distribution $F$. The main result for the model is a shape theorem
stating that, if

\begin{equation}\label{eq:rvillk}
\int_0^\infty e^{-\varphi r}dF(r)<\infty \textrm{ for some }
\varphi<0,
\end{equation}

\noindent then, on the scale $1/t$, the set $S_t$ has an
asymptotic shape, which by rotational invariance is forced to be a
Euclidean ball; see Theorem 1.1 in Deijfen \emph{et al} (2003).
This is a stronger result than for the Richardson model, for which
the existence of an asymptotic shape is established but the nature
of the shape so far remains unknown.\medskip

\noindent In H\"{a}ggstr\"{o}m and Pemantle (1998,2000) a
generalized version of the Richardson model is introduced that
contains two infection types, referred to as type 1 and type 2
infection respectively. The dynamics is that, for $i=1,2$, a
healthy site becomes type $i$ infected at rate $\lambda_i$ times
the number of type $i$ infected nearest neighbors and once
infected it never recovers. Here $\lambda_i\in\mathbb{R}$
indicates the intensity of the type $i$ infection. Furthermore,
when a site is struck by one of the infection types, it is assumed
to become immune to the other type, making the process a model for
competing growth on $\mathbb{Z}^d$.\medskip

\noindent The continuum model of Deijfen (2003) can also be
extended to incorporate two infection types; this is done in
Deijfen \emph{et al} (2003). Like in the one-type model, the
growth takes place by way of spherical outbursts, the type of an
outburst being determined by the region in which it takes place.
Let $S_t^1$ ($S_t^2$) denote the region occupied by type 1 (2)
infection at time $t$. The time until an outburst occurs in
$S_t^i$ ($i=1,2$) is exponentially distributed with parameter
$\lambda_i|S_t^i|$, where $\lambda_i\in \mathbb{R}$ is the
intensity of the type $i$ infection, and the outburst point is
chosen uniformly in $S_t^i$. An outburst in the type $i$ infected
region hurls out type $i$ infection to all previously uninfected
points in a ball centered at the outburst point. The radii in the
outburst balls are i.i.d.\ random variables with the same
distribution, denoted by $F$, for both outburst types. Note that
the infection types are assumed to be mutually exclusive, that is,
a point can be contained in at most one of the sets $S_t^1$ and
$S_t^2$.\medskip

\noindent According to the above mechanisms, the two-type models
describe a competition for space -- on $\mathbb{Z}^d$ and
$\mathbb{R}^d$ respectively -- between two infections. For both
models there are three possible results of this competition:
Either one of the infection types wins or there is a draw. For
$i=1,2$, define
$$
A_i=\{\textrm{the type $i$ infection reaches sites arbitrarily far
from the origin}\}.
$$
We say that the type 1 (2) infection has won if $\neg A_2$ ($\neg
A_1$) occurs and with a draw we mean that $A=A_1\cap A_2$ occurs,
that is, a power balance is established so that both infection
types keep growing indefinitely; here $\neg$ denotes set
complement.\medskip

\noindent In the discrete model and in the continuum model with
bounded outburst radii, clearly both events $\neg A_1$ and $\neg
A_2$ -- and thereby also the events $A_1$ and $A_2$ -- have
positive probability regardless of the values of the intensities
$\lambda_1$ and $\lambda_2$. In the continuum model with unbounded
support for the outburst radii, this fact requires a proof; see
Proposition 5.1 in Deijfen \emph{et al} (2003). The mathematical
challenge with the two-type models lies in deciding whether the
event $A$ has positive probability or not. Intuitively, the state
of equilibrium represented by $A$ should be possible to maintain
in the long run if and only if the infection types have the same
intensity. To state this conjecture more formally, let
$\tilde{P}^{\lambda_1,\lambda_2}$ denote the probability law for a
discrete two-type process started from two single nodes, one at
the origin and one at the point $(1,0,\ldots,0)$, and write
$P^{\lambda_1,\lambda_2}$ for the distribution of a continuum
two-type process started from two unit balls located next to each
other, one centered at the origin and the other at the point
$(2,0,\ldots,0)$.

\begin{conj}
For any dimension $d\geq 2$, we have
\begin{itemize}
\item[\rm{(a)}] $\tilde{P}^{\lambda_1,\lambda_2}(A)>0\Leftrightarrow
\lambda_1=\lambda_2$;
\item[\rm{(b)}] if $F$ satisfies (\ref{eq:rvillk}),
then $P^{\lambda_1,\lambda_2}(A)>0\Leftrightarrow
\lambda_1=\lambda_2$.
\end{itemize}
\end{conj}

\noindent \textbf{Remark 1.1} It is shown in Deijfen
and H\"{a}ggstr\"{o}m (2003) that the possibility of infinite
coexistence in the lattice model is not affected by the initial
configuration, that is, whether
$\tilde{P}^{\lambda_1,\lambda_2}(A)$ is positive or not does not
depend on the initial state of the model, given of course that
neither of the infection types surrounds the other already at the
beginning. In Deijfen \emph{et al} (2003) an analogous result
(stated as Proposition \ref{prop:startomr} in the present paper)
is proved for the continuum model provided the
technical condition
\begin{equation}  \label{eq:small_outbursts}
F(\varepsilon)>0 \mbox{ for all } \varepsilon>0 \, .
\end{equation}
Hence the particular choice of initial sets is unimportant in
Conjecture 1.1. $\hfill \Box$
\medskip

\noindent Efforts have been made to prove Conjecture 1.1 and a
number of partial results have been obtained. As for the case with
unequal intensities, it is shown in H\"{a}ggstr\"{o}m and Pemantle
(2000) that, if $\lambda_1$ is fixed, then
$\tilde{P}^{\lambda_1,\lambda_2}(A)=0$ for all but at most
countably many values of $\lambda_2$, and in Deijfen \emph{et al}
(2003) the same result is established for
$P^{\lambda_1,\lambda_2}$ under the assumption that
(\ref{eq:rvillk}) holds for $F$. Unfortunately, to strengthen
these results to cover \emph{all} intensities $\lambda_2\neq
\lambda_1$ seems to be marred by large technical
difficulties.\medskip

\noindent The discrete model with equal intensities was first
studied by H\"{a}ggstr\"{o}m and Pemantle (1998). Their main
result is that, if $\lambda_1=\lambda_2$ and $d=2$, then
$\tilde{P}^{\lambda_1,\lambda_2}(A)>0$, that is, on a
two-dimensional square lattice, mutual unbounded growth is possible for
infections with the same intensity. Recently Garet and Marchand
(2004) were able to extend this result to arbitrary dimensions.
While the technique of H\"{a}ggstr\"{o}m and Pemantle relies on some
very specific features of the square lattice $\mathbb{Z}^2$, the
proof of Garet and Marchand is much less sensitive to the exact
lattice structure, and thus more amenable to generalizations. In this
paper, we exploit the ideas of Garet and Marchand to establish infinite
coexistence in the continuum setting. That is, we prove the following.

\begin{theorem}\label{th:hres}
Let $d\geq 2$ and assume that $F$ satisfies conditions
(\ref{eq:rvillk}) and (\ref{eq:small_outbursts}).
Then, if $\lambda_1=\lambda_2$, we have
$P^{\lambda_1,\lambda_2}(A)>0$.
\end{theorem}

\medskip\noindent
\textbf{Remark 1.2} Condition (\ref{eq:small_outbursts}), which
says that arbitrarily small outbursts happen with positive
probability, is only needed in order to apply Proposition
\ref{prop:startomr} below (quoted from Deijfen \emph{et al}
(2003)). Without assuming (\ref{eq:small_outbursts}), we still
have the slightly weaker result that infinite coexistence has
positive probability for certain other compactly supported initial
states -- see Theorem \ref{th:hresmod}. Furthermore, we strongly
believe that Proposition \ref{prop:startomr} holds without the
condition (\ref{eq:small_outbursts}). $\hfill \Box$

\medskip\noindent
The most important step we need to take in order to be able to
invoke the Garet--Marchand ideas is a new coupling construction of
the two-type continuum model, designed in such a way that certain
subadditivity properties of the process are reflected in the
coupling. This construction is described in Section 2 together
with some preliminary results, and the proof of Theorem
\ref{th:hres} is given in Section 3. Superficially, our proof may
not seem to have much in common with the Garet--Marchand paper,
but the underlying key idea is the same. Translated to the
continuum setting, the basic idea is to observe (\ref{eq:limsup})
and to demonstrate that $P^{\lambda_1,\lambda_2}(A)=0$ would lead
to consequences that contradict (\ref{eq:limsup}).

\section{A new coupling construction}

\noindent Let us begin by describing in more detail how the
continuum two-type model with equally powerful infections can be
generated. To be able to exploit the Garet--Marchand ideas, we
need to modify the construction in Deijfen (2003) and Deijfen
\emph{et al} (2003) in order to make a certain subbadditivity
relation (\ref{eq:subadd}) hold.

\medskip\noindent
To construct the model in $d$ dimensions, consider a Poisson
process on $\mathbb{R}^{d+1}$ with rate $\lambda$, where the extra
dimension represents time and $\lambda$ denotes the common
intensity of the infection types. Associate independently to each
Poisson point a random radius variable with distribution $F$ and
pick two bounded disjoint regions $S_0^1,S_0^2\subset
\mathbb{R}^d$ to start up the growth, the set $S_0^i$ ($i=1,2$)
indicating the initially type $i$ infected region.\medskip

\noindent Starting at time 0, the growth is now brought about by
scanning within the region $(S_0^1\cup S_0^2)\times \mathbb{R}$
upwards along the time axis until a point in the Poisson process
is hit. An outburst then occurs at this point, causing the
previously uninfected parts of a ball around the outburst point to
be infected. The radius of the outburst ball is given by the
radius variable associated with the Poisson point and the type of
the transmitted infection is determined by in which one of the
regions $S_0^1\times \mathbb{R}$ and $S_0^2\times \mathbb{R}$ the
Poisson point is found. Let $T_1$ denote the time at which this
first outburst takes place and, for $i=1,2$, let $S_{(1)}^i$
denote the type $i$ infected region after the outburst. Here the
region with the same type as the outburst might be enlarged as
compared to before the outburst, while the other region remains
unchanged.

\medskip\noindent
To generate the next outburst, we keep scanning upwards in time
for Poisson points, but in such a way that the region $S_0^1\cup
S_0^2$ and the region $(S_{(1)}^1\cup S_{(1)}^2) \setminus
(S_0^1\cup S_0^2)$ are treated somewhat differently: for the
region $S_0^1\cup S_0^2$, we scan from time $T_1$ and upwards,
while simultaneously for $(S_{(1)}^1\cup S_{(1)}^2) \setminus
(S_0^1\cup S_0^2)$, we scan from time $0$ and upwards. Hence, the
time $T_2$ of the next outburst is given by the minimum of $T'_2$
and $T''_2$, where $T'_2$ is the time coordinate of the second
Poisson point in $(S_0^1\cup S_0^2)\times \mathbb{R}$, and $T''_2$
is the sum of $T_1$ and the time coordinate of the first Poisson
point in $(S_{(1)}^1\cup S_{(1)}^2) \setminus (S_0^1\cup
S_0^2)\times \mathbb{R}$.

\medskip\noindent
Subsequent outbursts are generated by similarly scanning up the
time axis. This means that a Poisson point at coordinates $(x, t)$
with $x \in \mathbb{R}^d$ and $t>0$ represents an outburst at $x$,
not at (absolute) time $t$, but precisely $t$ time units after the
point $x$ first became infected. We thus obtain two increasing
sequences $\{S^1_{(n)}\}$ and $\{S^2_{(n)}\}$ of regions in
$\mathbb{R}^d$ indicating the type $1$ and type $2$ infected
regions after $n$ outbursts, and a strictly increasing sequence
$\{T_n\}$ specifying the time points for the outbursts. For $t \in
[T_n, T_{n+1})$, the type $i$ infected region at time $t$ is given
by $S^i_t = S^i_{(n)}$ and the total infected region at time $t$
is $S_t = S^1_t \cup S^2_t$. Clearly this defines a model where
the time until an outburst occurs in $S_i^t$ is exponentially
distributed with parameter $\lambda|S^i_t|$ and the location of
the outburst is uniformly distributed in $S^i_t$, as
desired.\medskip

\noindent In summary, the probability space underlying our
construction, is simply a marked Poisson process on
$\mathbb{R}^{d+1}$ with rate $\lambda$ and i.i.d.\ marks with
distribution $F$. A more thorough description of the two-type
model can be found in Section 3 in Deijfen \emph{et al} (2003).
The construction there is different from the one given here, but
it is easy to see that, if the infection types have the same
intensity, the two construction methods result in growth processes
with the same distribution. The difference between the two
constructions matter only when we couple realizations of processes
starting from different initial configurations.\medskip

\noindent If we do not distinguish between the type 1 and the type
2 infection, then the symmetric two-type model reduces to the
one-type model of Deijfen (2003), that is, the total infected
region in the two-type model behaves like a one-type process. An
important consequence of this is the existence of a so called time
constant, denoted by $\mu$, indicating the inverse asymptotic
speed of the growth. Let $B(x,r)$ denote the closed ball with
radius $r$ centered at $x\in\mathbb{R}^d$, and, for $x,y\in
\mathbb{R}^d$, let $T_{x,y}$ denote the time when the ball
$B(y,1)$ is fully infected in a one-type process with rate
$\lambda$ started with $S_0=B(x,1)$. Also, write
$\mathbf{n}=(n,0,\ldots,0)$. If the radius distribution $F$
satisfies (\ref{eq:rvillk}), we have

\begin{equation}\label{eq:mu}
\lim_{n\rightarrow\infty}\frac{T_{\mathbf{0},\mathbf{n}}}{n}=\lambda^{-1}\mu\quad
\textrm{a.s. and in }L^1,
\end{equation}

\noindent and the same limit is obtained for the time when the
single point $\mathbf{n}$ is infected; see Deijfen \emph{et al}
(2003) or note that, with the modified construction of the
process, (\ref{eq:mu}) follows immediately from Liggett's
subadditive ergodic theorem. The reason for considering the time
when the entire unit ball around the target point is infected
(rather than the time at which the point itself is infected) is
that, in combination with the particular way the model is
constructed from the Poisson process, it gives rise to the
following very useful subadditivity property.

\begin{lemma}  \label{lem:subadditivity}
For any $x,y,z \in \mathbb{R}^d$, we have
\begin{equation}\label{eq:subadd}
T_{x,y}\leq T_{x,z}+T_{z,y}.
\end{equation}
\end{lemma}

\noindent \emph{Proof:} For $t\geq 0$, let $S_t$ denote the infected region
at time $t$ starting with $S_0=B(x,1)$, and let $\tilde{S}_t$ denote
the infected region at time $t$ starting with $\tilde{S}_0=B(z,1)$.
We need to show that

\begin{equation}  \label{eq:need_to_show}
B(y,1)\subseteq S_{T_{x,z}+T_{z,y}}.
\end{equation}

\noindent By definition of $T_{z,y}$, we have that
$B(y,1)\subseteq \tilde{S}_{T_{z,y}}$. Hence it suffices to show
that, for any $t \geq 0$, we have
\begin{equation}  \label{eq:set_inclusion}
\tilde{S}_t\subseteq S_{T_{x,z}+t} \, .
\end{equation}
For $t = 0$, (\ref{eq:set_inclusion}) is immediate from the
definition of $T_{x,z}$. It remains to show that the set inclusion
is preserved as $t$ increases. The only $t$'s for which it could
possibly stop holding are at the times of outbursts in the
$\tilde{S}$-process. Consider the first such outburst at which
(\ref{eq:set_inclusion}) ceases to hold, denote the time (in the
$\tilde{S}$-process) at which this outburst happens by $t'$, and
denote its location by $x$. Since (\ref{eq:set_inclusion}) holds
for all $t<t'$, we have that at time $T_{x,z}+t'$ the point $x$
has been infected in the $S$-process for at least as long as it
has been infected in the $\tilde{S}$-process at time $t'$. It
therefore follows from the construction of the process that the
$S$-process by time $T_{x,z}+t'$ has already had an outbreak at
the point $x$ (or has such an outbreak precisely at time
$T_{x,z}+t'$) of the same radius as the outbreak at $x$ in the
$\tilde{S}$-process. Hence, (\ref{eq:set_inclusion}) never ceases
to hold and consequently holds for all $t$. In particular, taking
$t=T_{z,y}$ yields (\ref{eq:need_to_show}), and we are done.
$\hfill \Box$

\medskip\noindent
Using the subadditivity in Lemma \ref{lem:subadditivity}, let
us prove the following lemma which will be crucial in the proof
of Theorem \ref{th:hres}.\medskip

\begin{lemma}\label{lemma} We have
\begin{equation}  \label{eq:limsup}
\limsup_{m\rightarrow \infty} {\rm{E}}[T_{\mathbf{n},\mathbf{-m}}-
T_{\mathbf{0},\mathbf{-m}}] \geq n\lambda^{-1}\mu.
\end{equation}
\end{lemma}

\noindent \emph{Proof:} Note first that, since
$\textrm{E}[T_{\mathbf{n},\mathbf{-m}}] =
\textrm{E}[T_{\mathbf{0},\mathbf{n+m}}]$ and
$\textrm{E}[T_{\mathbf{0},\mathbf{-m}}] =
\textrm{E}[T_{\mathbf{0},\mathbf{m}}]$, we have
\[
\textrm{E}[T_{\mathbf{n},\mathbf{-m}}- T_{\mathbf{0},\mathbf{-m}}]
= \textrm{E}[T_{\mathbf{0},\mathbf{n+m}}-
T_{\mathbf{0},\mathbf{m}}] \, .
\]
Hence we are done if we can show that
\begin{equation}\label{eq:limsup'}
\limsup_{m\rightarrow \infty}
\textrm{E}[T_{\mathbf{0},\mathbf{n+m}}- T_{\mathbf{0},\mathbf{m}}]
\geq n\lambda^{-1}\mu.
\end{equation}
To do this, note that, for any positive integer $k$,
$\textrm{E}[T_{\mathbf{0},k\mathbf{n}}]$ can be rewritten as the
telescoping sum
\begin{equation}  \label{eq:telescope}
\textrm{E}[T_{\mathbf{0},k\mathbf{n}}] = \sum_{i=1}^k
\textrm{E}[T_{\mathbf{0},i\mathbf{n}}-
T_{\mathbf{0},(i-1)\mathbf{n}}].
\end{equation}
Suppose now for contradiction that (\ref{eq:limsup'}) fails. Then
there exists an $\varepsilon > 0$ such that
$\textrm{E}[T_{\mathbf{0},i\mathbf{n}}-
T_{\mathbf{0},(i-1)\mathbf{n}}] < n\lambda^{-1}\mu - \varepsilon$
for all but at most finitely many positive integers $i$. Using
(\ref{eq:telescope}), this implies that
\begin{equation}  \label{eq:on_one_hand}
\limsup_{k\rightarrow \infty}
\frac{\textrm{E}[T_{\mathbf{0},k\mathbf{n}}]}{k} \leq
n\lambda^{-1}\mu - \varepsilon,
\end{equation}
but according to (\ref{eq:mu}) we have
\begin{equation}  \label{eq:and_on_the_other}
\lim_{k\rightarrow \infty}
\frac{\textrm{E}[T_{\mathbf{0},k\mathbf{n}}]}{k} =
n\lambda^{-1}\mu.
\end{equation}
Comparing (\ref{eq:on_one_hand}) and (\ref{eq:and_on_the_other})
yields the desired contradiction, and the lemma is proved. $\hfill
\Box$

\section{Proof of Theorem \ref{th:hres}}

\noindent First we quote a proposition that states that the
initial configuration is basically irrelevant for the possibility
of mutual unbounded growth in the two-type continuum model. Here
$P_{\Gamma_1,\Gamma_2}^{\lambda_1,\lambda_2}$ denotes the
probability law of a process started with $S^1_0=\Gamma_1$ and
$S^2_0=\Gamma_2$.

\begin{prop}[Deijfen \emph{et al} 2003]\label{prop:startomr}
Let $(\Gamma_1,\Gamma_2)$ and $(\Gamma_1',\Gamma_2')$ be two pairs
of disjoint, bounded subsets of $\mathbb{R}^d$ with strictly
positive Lebesgue measures. Furthermore, suppose that the radius
distribution $F$ has unbounded support and satisfies
$F(\varepsilon)>0$ for all $\varepsilon>0$. Then
$$
P_{\Gamma_1,\Gamma_2}^{\lambda_1,\lambda_2}(A)>0\Rightarrow
P_{\Gamma_1',\Gamma_2'}^{\lambda_1,\lambda_2}(A)>0.
$$
\end{prop}

\noindent \textbf{Remark 3.1} Proposition \ref{prop:startomr}
extends to the case with bounded support as well, provided the
following (obviously necessary) condition on $(\Gamma'_1,
\Gamma'_2)$: If the radius distribution is bounded by $r$, then
neither of $\Gamma'_1$ or $\Gamma'_2$ may contain an impenetrable
layer of thickness $r$ around the other.$\hfill \Box$\medskip

\noindent Now, to prove Theorem \ref{th:hres} we will show that
infinite coexistence is indeed possible for two infections with
the same intensity if the initial sets are located sufficiently
far away from each other; see Theorem \ref{th:hresmod}. In view of
the above proposition, clearly this implies Theorem \ref{th:hres}.
To simplify notation, write
$P_{B(\mathbf{0},1),B(\mathbf{n},1)}^{\lambda_1,\lambda_2}=
P_{\mathbf{0},\mathbf{n}}^{\lambda_1,\lambda_2}$.

\begin{theorem}\label{th:hresmod}
If $F$ satisfies (\ref{eq:rvillk}) and $\lambda_1=\lambda_2$,
then $P_{\mathbf{0},\mathbf{n}}^{\lambda_1,\lambda_2}(G)>0$ for
large $n$.
\end{theorem}

\noindent\emph{Proof:} By time-scaling, we may without loss of
generality assume that $\lambda_1=\lambda_2=1$. Hence, in what
follows, only unit rate processes are considered. To avoid
superfluous use of superscripts, write
$P_{\mathbf{0},\mathbf{n}}^{1,1}=P_{\mathbf{0},\mathbf{n}}$.\medskip

\noindent First, fix $\varepsilon>0$ and pick $n$ large so that

\begin{itemize}
\item[(i)] E$[T_{\mathbf{0},\mathbf{n}}]\leq (1+\varepsilon)n\mu$;
\item[(ii)] $P_{\mathbf{0},\mathbf{n}}(T_{\mathbf{0},\mathbf{n}}<
(1-\varepsilon)n\mu)<\varepsilon$,
\end{itemize}

\noindent which is possible because of (\ref{eq:mu}). The symbol
'E' will throughout the proof be used to denote expected value
with respect to the marked Poisson process underlying the coupling
construction described in Section 2. Now, by Lemma \ref{lemma},
there exist arbitrarily large $m$ such that

\begin{equation}\label{eq:vvdiffbeg}
\textrm{E}[T_{\mathbf{n},\mathbf{-m}}-T_{\mathbf{0},\mathbf{-m}}]\geq
(1-\varepsilon)n\mu.
\end{equation}

\noindent We will show that, if $P_{\mathbf{0},\mathbf{n}}(A)=0$
and $\varepsilon$ is small, then (\ref{eq:vvdiffbeg}) fails for
large $m$, implying that we must have
$P_{\mathbf{0},\mathbf{n}}(A)>0$.\medskip

\noindent Denote by $B_k$ a ball with radius $k$ centered at the
point $\mathbf{n}/2$ and, for a two-type process started with
$S_0^1=B(\mathbf{0},1)$ and $S_0^2=B(\mathbf{n},1)$, let
$$
A_k^i=\{\textrm{the type $i$ infection reaches points outside
$B_k$}\}.
$$
Also, define $A_k=A_k^1\cap A_k^2$, that is, $A_k$ is the event
that both infection types reach points outside $B_k$. Note that,
since $A=\cap_{k=1}^\infty A_k$ and $A_k\supset A_{k+1}$, we have
$$
P_{\mathbf{0},\mathbf{n}}(A)=\lim_{k\rightarrow\infty}P_{\mathbf{0},\mathbf{n}}(A_k).
$$
Now assume that $P_{\mathbf{0},\mathbf{n}}(A)=0$ and pick $m$
large enough so that $P_{\mathbf{0},\mathbf{n}}(A_m)<1/2$.
Trivially,

\begin{eqnarray*}
\textrm{E}[T_{\mathbf{n},\mathbf{-m}}-T_{\mathbf{0},\mathbf{-m}}]
& = & \textrm{E}\left[\mathbf{1}_{\{A_m^1\}}
\right(T_{\mathbf{n},\mathbf{-m}}-
T_{\mathbf{0},\mathbf{-m}}\left)\right]\\
& + & \textrm{E}\left[\mathbf{1}_{\{\neg
A_m^1\}}(T_{\mathbf{n},\mathbf{-m}}-
T_{\mathbf{0},\mathbf{-m}})\right],
\end{eqnarray*}

\noindent where $\mathbf{1}_{\{\cdot\}}$ denotes the indicator
function. On $\neg A_m^1$, we have
$T_{\mathbf{n},\mathbf{-m}}-T_{\mathbf{0},\mathbf{-m}}<0$, and
hence

\begin{eqnarray}  \nonumber
\textrm{E}[T_{\mathbf{n},\mathbf{-m}}-T_{\mathbf{0},\mathbf{-m}}]
& \leq &
\textrm{E}\left[\mathbf{1}_{\{A_m^1\}}(T_{\mathbf{n},\mathbf{-m}}-T_{\mathbf{0},\mathbf{-m}})\right]\\
& \leq & \textrm{E}\left[\mathbf{1}_{\{A_m^1\}}
T_{\mathbf{0},\mathbf{n}}\right], \label{eq:first_estimate}
\end{eqnarray}

\noindent where the last inequality follows from Lemma
\ref{lem:subadditivity}. To bound the right hand side from above,
write

\begin{equation}  \label{eq:rewritten_estimate}
\textrm{E}\left[\mathbf{1}_{\{A_m^1\}}T_{\mathbf{0},\mathbf{n}}\right]=
\textrm{E}[T_{\mathbf{0},\mathbf{n}}]-\textrm{E}\left[\mathbf{1}_{\{\neg
A_m^1\}}T_{\mathbf{0},\mathbf{n}}\right].
\end{equation}

\noindent By the choice of $n$, we have
$\textrm{E}[T_{\mathbf{0},\mathbf{n}}]\leq (1+\varepsilon)n\mu$,
and so it remains to find a lower bound for
$\textrm{E}[\mathbf{1}_{\{\neg
A_m^1\}}T_{\mathbf{0},\mathbf{n}}]$. To this end, remember that,
also by the choice of $n$, we have
$P_{\mathbf{0},\mathbf{n}}(T_{\mathbf{0},\mathbf{n}}<
(1-\varepsilon)n\mu)< \varepsilon$. Hence

\begin{eqnarray*}
\textrm{E}\left[\mathbf{1}_{\{\neg
A_m^1\}}T_{\mathbf{0},\mathbf{n}}\right] & \geq &
P_{\mathbf{0},\mathbf{n}}\left(\neg A_m^1\cap
\{T_{\mathbf{0},\mathbf{n}}
\geq (1-\varepsilon)n\mu\}\right)\cdot(1-\varepsilon)n\mu\\
& \geq & \left[P_{\mathbf{0},\mathbf{n}}(\neg
A_m^1)-\varepsilon\right]\cdot (1-\varepsilon)n\mu.
\end{eqnarray*}

\noindent By symmetry, $P_{\mathbf{0},\mathbf{n}}(\neg A_m^1)=
P_{\mathbf{0},\mathbf{n}}(\neg A_m^2)=
P_{\mathbf{0},\mathbf{n}}(\neg A_m)/2$, and, using the fact that
$P_{\mathbf{0},\mathbf{n}}(A_m)<1/2$, it follows that
$P_{\mathbf{0},\mathbf{n}}(\neg A_m^1)\geq 1/4$. Thus
$$
\textrm{E}\left[\mathbf{1}_{\{\neg
A_m^1\}}T_{\mathbf{0},\mathbf{n}}\right]\geq
\left(\frac{1}{4}-\varepsilon\right)(1-\varepsilon)n\mu,
$$
which in conjunction with
(\ref{eq:first_estimate}) and (\ref{eq:rewritten_estimate}) yields
the estimate

\begin{eqnarray*}
\textrm{E}[T_{\mathbf{n},\mathbf{-m}}-T_{\mathbf{0},\mathbf{-m}}]
& \leq & \textrm{E}[T_{\mathbf{0},\mathbf{n}}]
-\textrm{E}\left[\mathbf{1}_{\{\neg
A_m^1\}}T_{\mathbf{0},\mathbf{n}}\right] \\
& \leq & (1+\varepsilon)n\mu-\left(\frac{1}{4}-\varepsilon\right)(1-\varepsilon)n\mu\\
& \leq & \left(\frac{3}{4}+3\varepsilon\right)n\mu.
\end{eqnarray*}

\noindent Now, if $\varepsilon$ is small, clearly
$(3/4+3\varepsilon)<(1-2\varepsilon)$. Hence, to sum up, assuming
that $P_{\mathbf{0},\mathbf{n}}(A)=0$, we have showed that, if
$\varepsilon$ is small and $m$ is large, then
$\textrm{E}[T_{\mathbf{n},\mathbf{-m}}-T_{\mathbf{0},\mathbf{-m}}]<(1-2\varepsilon)n\mu$.
But this contradicts the fact that (\ref{eq:vvdiffbeg}) should
hold for arbitrarily large $m$ and consequently we must have
$P_{\mathbf{0},\mathbf{n}}(A)>0$. $\hfill \Box$

\section*{References}

\noindent Deijfen, M. (2003): Asymptotic shape in a continuum
growth model, \emph{Adv. Appl. Prob.} \textbf{35},
303-318.\medskip

\noindent Deijfen, M. and H\"{a}ggstr\"{o}m, O. (2003): The
initial configuration is irrelevant for the possibility of mutual
unbounded growth in the two-type Richardson model, \emph{Comb.
Prob. Computing}, to appear.\medskip

\noindent Deijfen, M., H\"{a}ggstr\"{o}m, O. and Bagley, J.
(2003): A stochastic model for competing growth on $\mathbb{R}^d$,
\emph{Markov Proc. Relat. Fields} \textbf{10:2}, 217-248.\medskip

\noindent Garet, O. and Marchand, R. (2004): Coexistence in
two-type first-passage percolation models, \emph{Ann. Appl.
Prob.}, to appear.\medskip

\noindent H\"{a}ggstr\"{o}m, O. and Pemantle, R. (1998): First
passage percolation and a model for competing spatial growth,
\emph{J. Appl. Prob.} \textbf{35}, 683-692.\medskip

\noindent H\"{a}ggstr\"{o}m, O. and Pemantle, R. (2000): Absence
of mutual unbounded growth for almost all parameter values in the
two-type Richardson model, \emph{Stoch. Proc. Appl.} \textbf{90},
207-222.\medskip
\end{document}